\documentclass[11pt]{article}

\usepackage{amscd,amsmath, amssymb, fancyhdr, epsfig, color, tocloft, mathtools}
\usepackage{dutchcal}

\usepackage[backref=page]{hyperref}
\renewcommand*{\backrefalt}[4]{%
	\ifcase #1 (Not cited.)%
	\or        (Cited on page~#2.)%
	\else      (Cited on pages~#2.)%
	\fi}

\hypersetup{
	colorlinks   = true,
	citecolor    = magenta,
	linkcolor    = blue,
	urlcolor     = magenta	
}

\newcommand{\version}{version 2.0,\ \ Aug. 24, 2022}

\makeatletter
\def\x@arrow{\DOTSB\Relbar}
\def\xlongequalsignfill@{\arrowfill@\x@arrow\Relbar\x@arrow}
\providecommand{\xlongequal}[2][]{%
	\ext@arrow 0099\xlongequalsignfill@{#1}{#2}}
\def\xlongrightarrowfill@{\arrowfill@\relbar\relbar\longrightarrow}
\makeatother

\numberwithin{equation}{section}

\newcommand{\Lie}{\operatorname{Lie}}

\def\eqref#1{(\ref{#1})}

\newcommand{\goth}{\mathfrak}

\newcommand{\C}{{\mathbb C}}
\newcommand{\R}{{\mathbb R}}

\def\1{\sqrt{-1}\,}

\newcommand{\cntrct}                
{\hspace{2pt}\raisebox{1pt}{\text{$\lrcorner$}}\hspace{2pt}}

\newcommand{\arrow}{{\:\longrightarrow\:}}

\renewcommand{\phi}{\varphi}
\renewcommand{\epsilon}{\varepsilon}
\renewcommand{\geq}{\geqslant}
\renewcommand{\leq}{\leqslant}

\newcommand{\const}{\operatorname{{\sf const}}}

\newcommand{\codim}{\operatorname{\sf codim}}
\renewcommand{\dim}{\operatorname{\sf dim}}

\newcommand{\kah}{{\operatorname{\text{\sf kah}}}}

\newcounter{Mycounter}[section]
\newcounter{lemma}[section]
\renewcommand{\thelemma}{{Lemma \thesection.\arabic{lemma}}}
\newcommand{\lemma}{%
	\setcounter{lemma}{\value{Mycounter}}
	\refstepcounter{lemma}
	\stepcounter{Mycounter}
	{\noindent \bf \thelemma:\ }}

\newcounter{claim}[section]
\newcounter{sublemma}[section]
\newcounter{corollary}[section]
\newcounter{theorem}[section]
\renewcommand{\thetheorem}{{Theorem \thesection.\arabic{theorem}}}
\newcommand{\theorem}{%
	\setcounter{theorem}{\value{Mycounter}}
	\refstepcounter{theorem}
	\stepcounter{Mycounter}
	{\noindent \bf \thetheorem:\ }}

\newcounter{conjecture}[section]
\newcounter{proposition}[section]
\renewcommand{\theproposition} {{Proposition \thesection.\arabic{proposition}}}
\newcommand{\proposition}{%
	\setcounter{proposition}{\value{Mycounter}}
	\refstepcounter{proposition}
	\stepcounter{Mycounter}
	{\noindent \bf \theproposition:\ }}

\newcounter{definition}[section]
\renewcommand{\thedefinition} {{Definition~\thesection.\arabic{definition}}}
\newcommand{\definition}{%
	\setcounter{definition}{\value{Mycounter}}
	\refstepcounter{definition}
	\stepcounter{Mycounter}
	{\noindent \bf \thedefinition:\ }}

\newcounter{example}[section]
\renewcommand{\theexample}{{Example \thesection.\arabic{example}}}
\newcommand{\example}{%
	\setcounter{example}{\value{Mycounter}}
	\refstepcounter{example}
	\stepcounter{Mycounter}
	{\noindent \bf \theexample:\ }}

\newcounter{remark}[section]
\renewcommand{\theremark}{{Remark \thesection.\arabic{remark}}}
\newcommand{\remark}{%
	\setcounter{remark}{\value{Mycounter}}
	\refstepcounter{remark}
	\stepcounter{Mycounter}
	{\noindent \bf \theremark:\ }}

\newcounter{problem}[section]
\newcounter{question}[section]
\makeatletter

\def\blacksquare{\hbox{\vrule width 5pt height 5pt depth 0pt}}
\def\endproof{\blacksquare}

\newcommand{\proof}{{\bf Proof: \ }}
\newcommand{\pstep}{{\bf Proof. Step 1: \ }}

\begin{document}

\begin{center}
{\Large\bf  A Calabi-Yau theorem for Vaisman manifolds}\\[5mm]
{\large
Liviu Ornea\footnote{Liviu Ornea is  partially supported by Romanian Ministry of Education and Research, Program PN-III, Project number PN-III-P4-ID-PCE-2020-0025, Contract  30/04.02.2021},  
Misha Verbitsky\footnote{Misha Verbitsky is partially supported by
 the HSE University Basic Research Program, FAPERJ E-26/202.912/2018 
and CNPq - Process 313608/2017-2.\\[1mm]

}\\[4mm]

}

\end{center}

{\small
\hspace{0.15\linewidth}
\begin{minipage}[t]{0.8\linewidth}
{\bf Abstract} \\ 
A compact complex Hermitian manifold $(M, I, \omega)$ 
is called { Vaisman} if $d\omega=\omega\wedge \theta$
and the 1-form $\theta$, called { the Lee form,} 
is parallel with respect to the Levi-Civita
connection. The volume form of $M$ is invariant
with respect to the action of the vector field
$X$ dual to $\theta$ (called { the Lee field})
and the vector field $I(X)$, called { the 
anti-Lee field}. The cohomology class of
$\theta$, called { the Lee class}, plays the
same role as the K\"ahler class in K\"ahler
geometry. We prove that a Vaisman metric
is uniquely determined by its volume
form and the Lee class, and, conversely,
for each Lee class $[\theta]$ and each 
Lee- and anti-Lee-invariant volume
form $V$, there exists a Vaisman structure
with the volume form $V$ and the 
Lee class $c[\theta]$. This is an analogue
of the Calabi-Yau theorem claiming that the
K\"ahler form is uniquely determined
by its volume and the cohomology class.
\end{minipage}
}

\tableofcontents

\section{Introduction}

E. Calabi
(\cite{_Calabi:ICM54_,_Calabi:vanishing_canonical_}) 
has noticed that a K\"ahler metric is uniquely
determined by its K\"ahler class and its volume form.
The Calabi conjecture, proven by S.-T. Yau some 23 years
later, claims that on any compact K\"ahler $n$-manifold $M$
there exists a unique K\"ahler metric $\omega$
with a given volume form $V$ and a given K\"ahler class
$[\omega]$ if $\int_M V= \int_M[\omega]^n$.\index{Calabi conjecture}
This statement is equivalent to the existence
and uniqueness of the solutions of the complex
Monge-Amp\`ere equation $(\omega +dd^c \phi)^n= V$.\index{Monge-Amp\`ere equation}

From this observation, Calabi obtained that (conditional
on the Calabi conjecture) every compact K\"ahler manifold
$M$ with $c_1(M)=0$ admits a unique Ricci-flat metric
in any given K\"ahler class; this Ricci-flat metric
is clearly Einstein. Finding and classifying the Einstein metrics 
(and, more generally, the ``extremal metrics'', defined\index{metric!extremal}
by Calabi in \cite{_Calabi:extremal_} as a generalization of Einstein
and constant scalar curvature K\"ahler metrics) 
was one of the central subjects of the complex
algebraic geometry since the 1980-ies.

The Calabi-Yau theorem has two facets, equally important.
One is related to finding the Ricci-flat metric on manifolds
with trivial canonical bundle. The other one is not
related to the canonical bundle in any way: it is 
a result which claims that a K\"ahler metric\index{theorem!Calabi-Yau}
is uniquely defined by its volume form and
the K\"ahler class. The existence of the Ricci-flat
metrics follows directly from this, more general, result.

We prove a version of the Calabi-Yau theorem for Vaisman
manifolds  (\ref{_CY_Vaisman_Theorem_}).

Let $M$ be a compact smooth
manifold, and $F\subset TM$ a smooth
foliation. It is called {\bf transversally K\"ahler}
if the normal bundle $TM/F$ is equipped with\index{foliation!transversally K\"ahler}
a Hermitian structure (that is, a complex structure
and a Hermitian metric) which is locally obtained
as the pullback of a K\"ahler structure on the leaf
space. Sasakian manifolds are prime examples of
transversally K\"ahler manifolds (the leaf space
of the Reeb foliation on a Sasakian manifolds
is K\"ahler).

A differential form on $M$ is called {\bf basic}
if it vanishes on $F$ and is locally obtained
as the pullback of a form on the leaf space.
The basic forms are preserved by de Rham differential,
and the cohomology of the basic forms is called
{\bf the basic cohomology}. \index{foliation!taut}

A foliation is {\bf taut} if the top basic cohomology
is non-zero (this is not the usual definition, but a theorem of
Habib and Richardson, see 
\cite{_Habib_Richardson:basic_}). 
For taut foliations, one has also
Poincar\'e duality on the basic cohomology, and the
identification between the basic cohomology and the basic
harmonic forms if a transversal Riemannian structure
is given.\index{cohomology!basic}

When $F$ is taut and transversally
K\"ahler, the basic cohomology should satisfy all the nice properties
of the cohomology of the K\"ahler manifolds: the $dd^c$-lemma,
the Hodge decomposition, the Hodge structure, Lefschetz
$\goth{sl}(2)$-action and so on. We proved it in the
situation when $F$ is trivialized by a group action
which preserves the transversally K\"ahler structure
(\cite[Theorem 5.4]{_ov_super_sas_}); 
a similar result was proven much earlier
by A. El Kacimi-Alaoui (\cite{_Kacimi_}).

In this situation, A. El Kacimi-Alaoui
proves the transversal Calabi-Yau theorem
showing that the transversal K\"ahler structures
are uniquely defined by the transversal volume
form and the transversal K\"ahler class. 
We give the uniqueness part of the proof
in \ref{_transversal_volume_defines_omega_0_Lemma_},
and for the existence refer to \cite{_Kacimi_}.

We apply these results to Vaisman geometry.
Any Vaisman manifold is equipped by 
a transversally K\"ahler, holomorphic 
foliation $\Sigma$ (\ref{_Subva_Vaisman_Theorem_}).
The transversally K\"ahler structure of this
foliation depends on the Vaisman metric,
however, the transversally complex structure
depends only on the complex structure of
the Vaisman manifold (\cite{tsuk} or 
\cite[Corollary 2.7]{tsu2}). We construct a
correspondence between the set of
Vaisman metrics on $M$ and the
set of transversal K\"ahler structures.

The transversal Calabi-Yau theorem implies 
an important result about the Vaisman metrics
(\ref{_CY_Vaisman_Theorem_}), showing that
the Vaisman metric is defined uniquely, up to a
constant multiplier, by the volume and the Lee class
$[\theta]\in H^1(M)$. The space of possible
Lee classes on $M$ (its ``Lee cone'')\index{Lee class}
is described in \cite{tsuk} (see also \cite{ov_lee}):
it is identified with a certain open half-space in $H^1(M,\R)$.
Then, similarly to the Calabi theorem parametrizing the K\"ahler
forms, the set of all Vaisman structures on $(M,I)$ is parametrized by
the cohomological data together with the volumes.

The Vaisman Calabi-Yau theorem is deduced
directly from the transversal Calabi-Yau theorem,
because the transversally K\"ahler form of
a Vaisman manifold, together with its Lee
class, uniquely defines the Vaisman structure
(\ref{_Vaisman_from_omega_0_and_Lee_Lemma_}).
On the other hand, the transversal volume form
uniquely defines, and is uniquely defined,
by the volume form of a Vaisman manifold 
\ref{_Vaisman_volume_from_omega_0_Lemma_}.
This is follows from a curious
observation, made by K. Tsukada
in \cite{tsu2}, who proved
that the direction of the Lee
field of a Vaisman manifold is
determined by its complex structure.

\section{Preliminaries}

Let $(M,I,g,\omega)$ be a Hermitian manifold, $\dim_\C M\geq 2$. Here $\omega(\cdot, \cdot)=g(I\cdot, \cdot)$.

\hfill

\definition
The Hermitian manifold $(M,I,g,\omega)$ is {\bf locally conformally K\"ahler} (LCK) if 
there exists a closed 1-form $\theta$ such that $d\omega=\theta\wedge\omega$.
The 1-form $\theta$ is called the {\bf Lee form} and the 
$g$-dual vector field $\theta^\sharp$ is called the {\bf Lee field}. 

\hfill

\remark \label{_Equiv_def_remark_}
One can easily see that this definition is equivalent with the existence of a K\"ahler cover $\Gamma\arrow(\tilde M,\tilde\omega)\arrow M$ such that the deck group $\Gamma$ acts by holomorphic homotheties (e. g. \cite{va_gd}).
Therefore, one can define a {\bf homothety character} $\chi:\Gamma \arrow \R^{>0}$ which associates to each deck transform $\gamma$ the scale factor $\frac{\gamma^*\tilde\omega}{\tilde \omega}$. When $\tilde M$ is the universal cover, the homothety character is a representation of $\pi_1(M)$ and uniquely defines the class $[\theta]\in H^1(M,\R)\simeq H_1(M)\simeq \frac{\pi_1(M)}{[\pi_1(M),\pi_1(M)]}$ of the Lee form.

\hfill

 In this note, we shall be interested in a particular
 subclass of LCK manifolds, namely the Vaisman manifolds. 
 
 \hfill
 
 \definition
 The LCK manifold $(M,I,g,\omega)$ is a {\bf Vaisman
   manifold} if the Lee form is parallel with respect to
 the Levi-Civita connection of the metric $g$.
 
 \hfill
 
 \example
 Almost all non-K\"ahler compact complex surfaces are LCK,
 see e.g. \cite{_ovv:surf_}. Diagonal Hopf surfaces and
 Kodaira surfaces are Vaisman, but Kato surfaces and Inoue
 surfaces are not Vaisman. In any dimension, all diagonal
 Hopf manifolds are Vaisman, see \cite{ov_pams}.
 
 \hfill
 
\remark
(\cite{va_gd})\label{_Theta_Killing_and_holo_Remark_} It
is easily seen that on a Vaisman manifold, the Lee and
anti-Lee vector fields $\theta^\sharp$ and
$I\theta^\sharp$ are Killing and holomorphic. Moreover,
they commute: $[\theta^\sharp,
  I\theta^\sharp]=0$. Therefore, they define a holomorphic
1-dimensional foliation $\Sigma$. One can also show that
$\Sigma$ is transversally K\"ahler and its leaves are
totally geodesic.

\hfill

\remark Up to a homothety, we can suppose that the length
of the Lee form is 1. Then (\cite[page 242]{va_gd}), the following identity holds:
\begin{equation}\label{_omega_via_theta_Chapter_8_Equation_}
	\omega= d^c \theta +\theta \wedge \theta^c.
\end{equation}

\hfill

\theorem \label{_Subva_Vaisman_Theorem_} 
Let $M$ be a compact Vaisman manifold, and 
$\Sigma\subset TM$ its canonical foliation. Then $\Sigma=\ker \omega_0$, where 
$\omega_0=d^c\theta$.

\proof 
\cite[Theorem 3.1]{va_gd} or \cite[Proposition 6.4]{_Verbitsky:Vanishing_LCHK_}. 
\endproof

\hfill

\remark (\cite{va_gd}) 
The K\"ahler cover of a compact Sasakian manifold is a
cone $S\times \R^{>0}$, where $S$ is a Sasakian manifold,
with cone metric $\tilde g= t^2g^S+dt\otimes dt$, where
$g^S$ is the Sasaki metric on $S$ and $t$ is the
coordinate on $ \R^{>0}$.
 
\hfill

In LCK geometry, the analogue of the K\"ahler cone of a
compact K\"ahler manifold is the ``Lee cone'', i.e. the
set of  classes in $H^1(M)$ which can be represented by
Lee forms of an LCK structure on the fixed complex
manifold $(M,I)$ admitting LCK structures. For Vaisman
manifolds, the ``Lee cone'' is known:

\hfill

\theorem (\cite{tsuk, ov_lee})\label{_Lee_cone_on_Vaisman_Theorem_}
Let $M$ be a compact Vaisman manifold. Then:
\begin{description}
\item[(i)] \[ H^1(M)= H_d^{1,0}(M) \oplus \overline{H_d^{1,0}(M)} \oplus \langle \theta\rangle\]
\item[(ii)] Consider a 1-form $\mu\in H^1(M)^*$ vanishing on 
$H_d^{1,0}(M) \oplus \overline{H_d^{1,0}(M)} \subset H^1(M)$
and satisfying $\mu([\theta])>0$, where $H_d^{1,0}(M)$
is the space of closed holomorphic 1-forms. Then 
a class $\alpha\in H^1(M,\R)$ is a 
Lee class for some LCK structure if and only if $\mu(x) >0$.
\end{description}
\endproof

\hfill

We shall need the following description of basic cohomology on compact Vaisman manifolds:

\hfill

\theorem 
(\cite{va_gd}, \cite{_ov_super_sas_})\label{_Vaisman_harmonic_forms_Theorem_}
Let $(M,I,g,\theta)$ be a compact Vaisman manifold of complex dimension $n$, with fundamental form $\omega$, and canonical foliation $\Sigma$. 
Denote by ${\cal H}^i$ the space of all  basic $i$-forms 
$\alpha\in \Lambda^*_\kah(M)$
which satisfy:
\begin{description}
	\item[\hspace{.1in}for $i \leq n$:]
	$\alpha$ is basic harmonic (i.\,e.    $\Delta_\kah(\alpha)=0$) and satisfies
	$\Lambda_{\omega_0}(\alpha)=0$; 
	\item[\hspace{.1in}for $i > n$:]
	$\alpha = \beta \wedge I\theta$
	where $\beta$ is basic harmonic and satisfies
	$L_{\omega_0}(\beta)=0$. 
\end{description}
Then all elements of ${\cal H}^*\oplus \theta\wedge {\cal H}^*$ are harmonic and, moreover,
all harmonic forms on $M$ belong to ${\cal H}^*\oplus \theta\wedge {\cal H}^*$.
\endproof

\section{The Lee field on a compact 
	Vaisman manifold} 

In \cite{tsu2}, K. Tsukada has shown that the
direction of the Lee field is uniquely determined by the
complex structure of the Vaisman manifold. We give
another proof of this result here; for other proofs,
see \cite{_Madani_Moroianu_Pilca:holo_} and 
\cite{_OV:EW_}.\footnote{The proof in \cite{_OV:EW_}
	is valid only for $b_1(M)=1$.}

\hfill

\proposition\label{_Lee_field_direction_Proposition_}
Let $M$ be a compact complex manifold of Vaisman type,
and $\theta^\sharp$ the Lee field of a Vaisman structure
$(\omega, \theta)$. Then $\theta^\sharp$ is determined by the
complex structure on $M$ uniquely up to a real multiplier.

\hfill

\pstep
Consider the form $\omega_0 = d\theta^c$.
As shown in \cite{_Verbitsky:Vanishing_LCHK_},
$\omega_0$ is a semi-positive Hermitian form, and
its kernel is precisely 
$\Sigma= \langle\theta^\sharp,I(\theta^\sharp)\rangle$.
Let $\omega_0'$ be a form associated in the same way 
with some other Vaisman structure $(\omega', \theta')$ on $M$.
Then $\eta:=\omega_0 + \omega_0'$ is an exact,
semi-positive (1,1)-form. This form cannot
be strictly positive because $\int_M \eta^{\dim_\C M}=0$,
hence it has a non-trivial kernel, which is
contained in $\ker \omega_0 \cap \ker \omega_0'$.
However, $\dim_\C \ker \omega_0 =\dim_\C\ker \omega_0'=1$,
and therefore these kernel spaces coincide.
This implies that the canonical foliation $\Sigma'$
associated with $(\omega', \theta')$ is equal to $\Sigma$.

\hfill

{\bf Step 2:}
Recall that $\theta^\sharp$ is 
holomorphic and Killing by \ref{_Theta_Killing_and_holo_Remark_}.
Since $\Sigma$ has a non-degenerate holomorphic section $\theta^\sharp$, it is
trivial as a holomorphic line bundle. Then $H^0(M, \Sigma)=\C$,
and the space of real holomorphic vector fields tangent
to $\Sigma$ has real dimension 2. Since $\theta^\sharp$
is Killing, it acts conformally on the 
K\"ahler covering $\tilde M$ of $M$.
However, a holomorphic conformal vector field
on a K\"ahler manifold multiplies the K\"ahler
form by a constant.
This gives a character $\sigma:\; \goth s \arrow \R$
on the Lie algebra of all holomorphic 
vector fields tangent to $\Sigma$.
We claim that $\sigma$ is uniquely
determined by the cohomology class of $\theta$.
Indeed, let $f\in C^\infty\tilde M$ be a function
such that $df =\theta$, and $X\in \goth s$, 
$X= a\theta^\sharp + bI(\theta^\sharp)$. Then 
\[ \Lie_X f= \langle \theta, X\rangle = 
a  \langle \theta, \theta^\sharp\rangle = a,\]
because $\nabla \theta=0$, and $a$
is equal to $\sigma(X)$ since
$\sigma(X)\tilde \omega =\Lie_X\tilde \omega = 
a \Lie_{\theta^\sharp}\tilde \omega+ 
b \Lie_{I(\theta^\sharp)}\tilde \omega=a\tilde \omega$.

\hfill

{\bf Step 3:} 
The anti-Lee field $I(\theta^\sharp)$
is distinguished by $\sigma(I(\theta^\sharp))=0$,
and the direction of $\theta^\sharp$ is determined by
$\sigma(\theta^\sharp)>0$. Therefore, to show that
$\theta^\sharp$ is independent from the choice
of the Vaisman structure, it would suffice to
prove that $\sigma$ is independent.

Let $\theta_1, \theta_2$ be two Lee classes of the
Vaisman metrics $\omega_1, \omega_2$ on $M$. 
Using the harmonic decomposition for 1-forms
(\cite{tsuk}), 
we obtain that any harmonic 1-form on $(M, \omega_1)$
is proportional to $\theta_1 + \alpha$, where
$\alpha$ is $\Sigma$-basic. The set
of possible Lee classes for Vaisman structures 
is a half-space in $H^1(M,\R)$ (\ref{_Lee_cone_on_Vaisman_Theorem_}),
with the boundary represented by $\Sigma$-basic forms.
Therefore, we can always replace the Lee class $[\theta]$ of a Vaisman
manifold by $[\const \theta]$, where $\const$ is positive.
Then, for some positive real constant $A$, the class
$[A\theta_1-\theta_2]$ is equal to $[\alpha]$, 
where $\alpha$ is a $\Sigma$-basic closed 1-form. 
Since the form $\alpha$ is basic, 
$\langle X, \alpha\rangle=0$ for any
$X$ tangent to $\Sigma$. Therefore, 
the character $\sigma_1:\; \goth s \arrow \R$ associated to
$[\theta_1]$ is proportional to the character
$\sigma_2:\; \goth s \arrow \R$  associated to
$[\theta_2]=[A\theta_1+\alpha]$.
\endproof


\section{The complex Monge-Amp\`ere equation}


We start by introducing a version of a result of
\cite{_OV:EW_} which proves the uniqueness of the solution
of the Monge-Amp\`ere equation on Vaisman manifolds.
Using a theorem of A. El Kacimi-Alaoui (\cite{_Kacimi_}),
we show that the solution always exists. This result
is a complete analogue of the existence and uniqueness
of the solutions of the complex Monge-Amp\`ere equations
on a compact K\"ahler manifold proven by S.-T. Yau.\index{Monge-Amp\`ere equation}

\hfill

We say that a form $\eta$ on a Vaisman manifold
{\bf is Lee-invariant} if  $\Lie_{\theta^{\sharp}}\eta=0$
and {\bf anti-Lee invariant} if
$\Lie_{I\theta^{\sharp}}\eta=0$.
The following theorem claims that the Vaisman
structure is uniquely determined by the
cohomological data and 
the Lee- and anti-Lee-invariant volume form.

\hfill

\theorem\label{_CY_Vaisman_Theorem_}
Let $(M,\omega,\theta)$ be a compact Vaisman manifold,
and $V'$ a Lee- and anti-Lee-invariant volume form on $M$,
satisfying $\int_M V' =\int_M \omega^n$.
Then there exists a unique Vaisman metric $\omega'$ 
on $M$ with the same Lee class
and the volume form $(\omega')^n=V'$.

\hfill

This theorem 
is proven later in this section.

\hfill

Recall that 
{\bf a transversal volume form}, or {\bf basic volume form}
is a basic form $V\in \Lambda^k_F(M)$ on a foliated manifold 
$(M, F)$, $k= \codim F$, which defines a non-degenerate volume
form locally on the leaf spaces of $F$. The following lemma
is a transversal form of Calabi-Yau theorem, essentially
due to \cite[\S 3.5.5]{_Kacimi_}.

\hfill

\lemma\label{_transversal_volume_defines_omega_0_Lemma_}
Let $M$ be a compact Vaisman $n$-manifold,
and $\Sigma$ its canonical foliation. 
Then for any $\Sigma$-basic volume form $V_0$ which is
cohomologous to an $(n-1)$-th power
of a transversally K\"ahler form $\eta_1$, there
exists a unique transversally K\"ahler form $\eta_2$
in the same basic cohomology class 
such that $\eta_2^{n-1}=V_0$.

\hfill

\proof
We start by proving the uniqueness of a transversally
K\"ahler form with a given transversal volume.
Using the transversal $dd^c$-lemma (\cite{_Kacimi_}),
we obtain
$\eta_1=\eta_2 + dd^c f$, where $f$ is a $\Sigma$-basic
function (that is, a function which is constant on the leaves of $\Sigma$).
Then $\eta_1^{n-1}-\eta_2^{n-1}= dd^c f\wedge P$, where
$P=\sum_{i=0}^{n-2} \eta_1^{i}\wedge \eta_2^{n-2-i}$.
Consider the operator
\[
f \mapsto {D_P}(f):=\frac{dd^c f\wedge P}{\eta_1^{n-1}}
\]
taking basic functions to basic functions. 
Let $U\subset M$ be a sufficiently small open set, and $X_U$
the leaf space of $\Sigma$ on $U$.
Clearly, the map $D_P:\; C^\infty(X_U) \arrow C^\infty(X_U)$ 
is a second order elliptic operator. By Hopf maximum
principle, 
 any non-constant $f \in \ker D_P$ cannot have a maximum.
However, any $\Sigma$-basic function on $M$ has a maximum
somewhere, because $M$ is compact. Therefore,
any $f\in \ker D_P$ is constant. This proves the
uniqueness of solutions. The existence of solutions
is obtained by repeating Yau's argument in the transversal setup,
as done in \cite[\S 3.5.5 (iv)]{_Kacimi_}.
\endproof

\hfill

\lemma\label{_Vaisman_from_omega_0_and_Lee_Lemma_}
Let $M$ be a compact complex $n$-manifold of Vaisman type.
Then a Vaisman structure on $M$ is uniquely determined
by its transversal K\"ahler form $\omega_0$ (\ref{_Subva_Vaisman_Theorem_})
and the Lee class $[\theta]\in H^1(M, \R)$.\footnote{The Lee
	class $[\theta]$ is uniquely determined by the homothety character $\chi$
	(\ref{_Equiv_def_remark_}) 
	and determines it.}

\hfill

\proof
By \eqref{_omega_via_theta_Chapter_8_Equation_}, 
we have $\omega=\omega_0 + \theta \wedge\theta^c$. 
Therefore, it would suffice
to show that the Lee form $\theta$ is uniquely determined 
by $\omega_0$ and the Lee class.
Let $\theta$ and $\theta'$ be two Lee forms of Vaisman manifolds,
with the same transversal K\"ahler form $\omega_0$. 
Denote by $\eta$ the 1-form $\theta- \theta'$.
Since $\omega_0 = d^c \theta=  d^c \theta'$ 
this would imply $d^c\eta= d\eta=0$. 
Such a 1-form cannot be exact, because 
if $\eta=df$, one has $dd^c f=0$;\index{maximum principle}
however, pluriharmonic functions are constant on any
compact manifold by the maximum principle. Therefore, $\theta$ 
cannot be cohomologous to $\theta'$.\index{function!pluriharmonic}
\endproof

\hfill

\lemma\label{_Vaisman_volume_from_omega_0_Lemma_}
Let $(M, \omega, \theta)$ be a Vaisman $n$-manifold,
and $\omega_0$ its transversal K\"ahler form (\ref{_Subva_Vaisman_Theorem_}).
Then $i_{I(\theta^\sharp)} i_{\theta^{\sharp}} \omega^n = n \omega_0^{n-1}$,
where $\theta^\sharp$ is the Lee field.

\hfill

\proof
By \eqref{_omega_via_theta_Chapter_8_Equation_}, 
we have $\omega=\omega_0 + \theta \wedge\theta^c$. 
This implies $\omega^n = \omega_0^{n-1} \wedge \theta \wedge\theta^c$.
The form $\omega_0$ is $\Sigma$-basic, and $|\theta^{\sharp}|=1$.
Then 
\[ 
i_{I(\theta^\sharp)} i_{\theta^{\sharp}} \omega^n = 
i_{I(\theta^\sharp)}i_{\theta^{\sharp}} (\theta \wedge\theta^c)\wedge \omega_0^{n-1}
=n\omega_0^{n-1}.\ \ \  \endproof
\]

Now we can prove \ref{_CY_Vaisman_Theorem_}.

The transversally K\"ahler form $\omega_0$ 
associated with the Vaisman structure
is uniquely defined by its
transversal volume form $V_0$ and its basic cohomology class
$[\omega_0]\in H^2_\Sigma(M)$, assuming that the transversal
cohomology classes of $\omega^{n-1}$ and $V$ are equal
(\ref{_transversal_volume_defines_omega_0_Lemma_}). However,
the class $[\omega_0]$ generates the kernel of the natural map
$H^2_\Sigma(M)\arrow H^2(M)$ (\ref{_Vaisman_harmonic_forms_Theorem_}).
This determines the basic cohomology class of $\omega_0$ up to a constant;
the constant is fixed if the transversal volume is fixed.
However, the transversal volume form $\omega_0^{n-1}$
is determined uniquely, up to a constant multiplier,
by the volume form 
$V$, as follows from \ref{_Vaisman_volume_from_omega_0_Lemma_}.
We obtain that the volume form
of a Vaisman metric uniquely defines $\omega_0$, in such a way
that $\omega_0^{n-1}=\const i_{I(\theta^\sharp)} i_{\theta^{\sharp}} V$
and the Lee class of the Vaisman structure associated with $\omega_0$
is proportional to $[\theta]$.

By \ref{_Lee_field_direction_Proposition_}, 
the Lee field $\theta^\sharp$ of a Vaisman structure is uniquely
(up to a constant) determined by the complex structure of $M$.
By \ref{_Vaisman_volume_from_omega_0_Lemma_}, 
$V_0= n^{-1} i_{I(\theta^\sharp)}i_{\theta^\sharp}V$ uniquely (up to a constant) defines
the transversal volume form  $\omega_0^{n-1}$, hence $\omega_0$
exists and is uniquely defined, up to a
constant.\footnote{This is where we use the assumption that $V$ is Lee- and
  anti-Lee-invariant; this is equivalent to $V_0$ being basic.}
Now, the Vaisman metric $\omega= \omega_0+ \theta \wedge \theta^c$
is uniquely defined by
$\omega_0$ and the Lee class (\ref{_Vaisman_from_omega_0_and_Lee_Lemma_}),
hence the Vaisman metric $\omega$ is uniquely 
(up to a constant) defined by its volume form 
$V$ and the Lee class. The constant is also fixed,
because $V= \omega^n$.

This proves uniqueness of a Vaisman metric with prescribed volume.
To see that a metric with a prescribed volume
form $V'$ exists, we write the corresponding transversal
volume form $V_0':= n^{-1} i_{I(\theta^\sharp)}i_{\theta^\sharp}V'$
and solve the transversal Calabi-Yau equation 
(\ref{_transversal_volume_defines_omega_0_Lemma_}), arriving at
a transversal K\"ahler form $\omega_0'= d^c\theta+ dd^c f =d^c\theta'$ satisfying
$(\omega_0')^{n-1}=V_0'$. Then 
$\omega':=\omega_0' + \theta'\wedge (\theta')^c$ is an LCK form 
which is invariant under the Lee field action; by \cite[Theorem A]{kor},
any LCK metric $\omega'$ admitting a conformal holomorphic
flow, non-isometric on its K\"ahler covering, is Vaisman.
\endproof

\hfill

{\bf Acknowledgements:} We are grateful to Nikita Klemyatin for 
finding errors in an earlier version.

\hfill

{\scriptsize

\noindent {\sc Liviu Ornea\\
University of Bucharest, Faculty of Mathematics and Informatics, \\14
Academiei str., 70109 Bucharest, Romania}, and:\\
{\sc Institute of Mathematics ``Simion Stoilow" of the Romanian
Academy,\\
21, Calea Grivitei Str.
010702-Bucharest, Romania\\
\tt lornea@fmi.unibuc.ro,   liviu.ornea@imar.ro}

\hfill

\noindent {\sc Misha Verbitsky\\
Instituto Nacional de Matem\'atica Pura e
Aplicada (IMPA) \\ Estrada Dona Castorina, 110\\
Jardim Bot\^anico, CEP 22460-320\\
Rio de Janeiro, RJ - Brasil }, also:\\
{\sc Laboratory of Algebraic Geometry, \\
Faculty of Mathematics, National Research University 
Higher School of Economics,
6 Usacheva Str. Moscow, Russia\\
\tt verbit@verbit.ru, verbit@impa.br }}


{\scriptsize
\begin{thebibliography}{199}
	
\bibitem[Ca1]{_Calabi:ICM54_} E. Calabi, {\em The Space of K\"ahler Metrics}, Proc. of the Int. Congress of Mathematicians 1954, Volume II,  206-207, E.P. Noordhoff, Groningen, 1956.
	
\bibitem[Ca2]{_Calabi:vanishing_canonical_} E. Calabi, 
	{\em On K\"ahler manifolds with vanishing canonical class} in Algebraic geometry and topology. A symposium in honor of S. Lefschetz,  78-89. Princeton Univ. Press, Princeton, N.J., 1957

\bibitem[Ca3]{_Calabi:extremal_} E. Calabi, {\em Extremal K\"ahler metrics}, in: Seminar on Differential Geometry, ed. S. T. Yau, Annals of Math. Studies 102, Princeton Univ. Press, Princeton, NJ (1982), 259-290.

\bibitem[El]{_Kacimi_}  A. El Kacimi-Alaoui, {\em Op\'erateurs transversalement elliptiques sur un feuilletage riemannien et applications}, Compositio Math. {\bf 73} (1990), no. 1, 57-106.

\bibitem[HR]{_Habib_Richardson:basic_} G. Habib, K. Richardson, {\em Modified differentials and basic cohomology for Riemannian foliations}, J. Geom. Anal. 23 (2013), no. 3, 1314-1342.

\bibitem[KO]{kor} 
Y. Kamishima, L. Ornea, {\em Geometric flow on compact  locally conformally K\"ahler manifolds}, Tohoku Math. J., {\bf 57} (2) (2005), 201-221.


\bibitem[MMP]{_Madani_Moroianu_Pilca:holo_} F. Madani, A. Moroianu, M. Pilca, {\em LCK structures with holomorphic lee vector field on
	vaisman-type manifolds}, Geom. Dedicata {\bf 213} (2021), 251-266. 
arXiv:1905.07300v1.  

\bibitem[OV1]{_OV:EW_} L. Ornea, M. Verbitsky, {\em Einstein-Weyl structures on complex manifolds and  conformal version of Monge-Amp\`ere equation}, Bull. Math. Soc. Sci. Math. Roumanie {\bf 51 (99)} No. 4, (2008), 339-353.

\bibitem[OV2]{ov_pams} L. Ornea, M. Verbitsky, {\em Locally conformally K\"ahler metrics obtained from pseudoconvex shells}, Proc. Amer. Math. Soc. {\bf 144} (2016), 325-335.

\bibitem[OV3]{_ov_super_sas_} 
 L. Ornea, M. Verbitsky,
  {\em Supersymmetry and Hodge theory on Sasakian and
    Vaisman manifolds}, arXiv:1910.01621, Manuscripta
  Math. \url{https://doi.org/10.1007/s00229-021-01358-8}


\bibitem[OV4]{ov_lee}  L. Ornea, M. Verbitsky, {\em Lee classes on LCK manifolds with potential}, arXiv:2112.03363.

\bibitem[Ts1]{tsuk} K. Tsukada, {\em Holomorphic forms and holomorphic vector fields on compact generalized Hopf manifolds}, Compositio Math. {\bf 93} (1994), no. 1, 1-22.

\bibitem[Ts2]{tsu2} 
K. Tsukada, 
{\em Holomorphic maps of compact generalized Hopf manifold
s}, Geom. Dedicata {\bf 68} (1997), 61-71.


\bibitem[Va]{va_gd} I. Vaisman, {\em Generalized Hopf manifolds}, Geom. Dedicata, {\bf 13} (1982), 231-255.

\bibitem[Ve]{_Verbitsky:Vanishing_LCHK_} 
M. Verbitsky, {\em Theorems on the vanishing of cohomology
for locally conformally hyper-K\"ahler manifolds},
Proc. Steklov Inst. Math. no. 3 ({\bf 246}), 54-78 (2004).


\bibitem[VVO]{_ovv:surf_} 
M. Verbitsky, V. Vuletescu, 
L. Ornea {\em Classification of non-K\"ahler surfaces and 
	locally conformally K\"ahler geometry}, Russian Math. Surv. {\bf 76} (2021), 261-290. arxiv:1810.05768.
	


\end{thebibliography}

}
\end{document}